# TAIL OF A LINEAR DIFFUSION WITH MARKOV SWITCHING

By Benoîte de Saporta and Jian-Feng Yao

*Université de Rennes 1*

Let $Y$ be an Ornstein–Uhlenbeck diffusion governed by a stationary and ergodic Markov jump process $X$: $dY_t = a(X_t)Y_t \, dt + \sigma(X_t) \, dW_t$, $Y_0 = y_0$. Ergodicity conditions for $Y$ have been obtained. Here we investigate the tail propriety of the stationary distribution of this model. A characterization of either heavy or light tail case is established. The method is based on a renewal theorem for systems of equations with distributions on $\mathbb{R}$.

**1. Introduction.** The discrete-time models $Y = (Y_n, n \in \mathbf{N})$ governed by a switching process $X = (X_n, n \in \mathbf{N})$ fit well to the situations where an autonomous process $X$ is responsible for the dynamic (or *regime*) of $Y$. These models are parsimonious with regard to the number of parameters, and extend significantly the case of a single regime. Among them, the so-called Markov-switching ARMA models are popular in several application fields, for example, in econometric modeling [see Hamilton (1989, 1990)]. More recently, continuous-time versions of Markov-switching models have been proposed in Basak, Bisi and Ghosh (1996) and Guyon, Iovleff and Yao (2004), where ergodicity conditions are established. In this paper we investigate the tail property of the stationary distribution of this continuous-time process. One of the main results (Theorem 2) states that this model can provide heavy tails, which is one of the major features required in nonlinear time-series modeling. Note that heavy tails may also be obtained by using a Lévy-driven Ornstein–Uhlenbeck (O.U.) process (without Markov switching); see Barndorff-Nielsen and Shephard (2001) and Brockwell (2001).

The considered process $Y$, called *diffusion with Markov switching*, is constructed in two steps:

---

---

First, the *switching process* $X = (X_t)_{t \geq 0}$ is a Markov jump process [see Feller (1966)], defined on a probability space $(\Omega, \mathcal{A}, Q)$, with a finite state space $E = \{1, \ldots, N\}$, $N > 1$. We assume that the intensity function $\lambda$ of $X$ is positive and the jump kernel $q(i, j)$ on $E$ is irreducible and satisfies $q(i, i) = 0$, for each $i \in E$. The process $X$ is ergodic and will be taken stationary with an invariant probability measure denoted by $\mu$.

Second, let $W = (W_t)_{t \geq 0}$ be a standard Brownian motion defined on a probability space $(\Theta, \mathcal{B}, Q')$, and let $\mathcal{F} = (\mathcal{F}_t)$ be the filtration of the motion. We will consider the product space $(\Omega \times \Theta, \mathcal{A} \times \mathcal{B}, (Q_x \otimes Q'))$, $\mathbb{P} = Q \otimes Q'$ and $\mathbb{E}$ the associated expectation. Conditionally to $X$, $Y = (Y_t)_{t \geq 0}$ is a real-valued diffusion process, defined, for each $\omega \in \Omega$ by:

1. $Y_0$ is a random variable defined on $(\Theta, \mathcal{B}, Q')$, $\mathcal{F}_0$-measurable;
2. $Y$ is solution of the linear SDE

$$\text{(1)} \qquad dY_t = a(X_t)Y_t \, dt + \sigma(X_t) \, dW_t, \qquad t \geq 0.$$

Thus $(Y_t)$ is a linear diffusion driven by an "exogenous" jump process $(X_t)$.

We say a continuous- or discrete-time process $S = (S_t)_{t \geq 0}$ is *ergodic* if there exists a probability measure $m$ such that when $t \to \infty$, the law of $S_t$ converges weakly to $m$ independently of the initial condition $S_0$. The distribution $m$ is then the *limit law* of $S$. When $S$ is a Markov process, $m$ is its unique invariant law.

In Guyon, Iovleff and Yao (2004), it is proved that the Markov-switching diffusion $Y$ is ergodic under the condition

$$\text{(2)} \qquad \alpha = \sum_{i \in E} a(i)\mu(i) < 0.$$

The main results of the present paper are the following theorems. Note that Condition 2 will be assumed satisfied throughout the paper and we denote by $\nu$ the stationary (or limit) distribution of $Y$.

THEOREM 1 (Light tail case).   *If for all $i$, $a(i) \leq 0$, then the stationary distribution $\nu$ of the process $Y$ has moments of all order; that is, for all $s > 0$ we have*

$$\int_{\mathbb{R}} |x|^s \nu(dx) < \infty.$$

THEOREM 2 (Heavy tail case).   *If there is an $i$ such that $a(i) > 0$, one can find an exponent $\kappa > 0$ and a constant $L > 0$ such that the stationary distribution $\nu$ of the process $Y$ satisfies*

$$t^\kappa \nu(]t, +\infty[) \underset{t \to +\infty}{\longrightarrow} L,$$

$$t^\kappa \nu(]-\infty, -t[) \underset{t \to +\infty}{\longrightarrow} L.$$

Note that the two situations from Theorems 1 and 2 form a dichotomy. Moreover, the characteristic exponent $\kappa$ in the heavy tail case is completely determined as following. Let

$$s_1 = \min\left\{\frac{\lambda(i)}{a(i)}\Big| a(i) > 0\right\},$$

$$M_s = \left(q(i,j)\frac{\lambda(i)}{\lambda(i) - sa(i)}\right)_{i,j \in E} \qquad \text{for } 0 \le s < s_1.$$

Then $\kappa$ is the unique $s \in ]0, s_1[$ such that the spectral radius of $M_s$ equals to 1.

The proof of Theorem 1 is a consequence of a result of Guyon, Iovleff and Yao (2004), and the proof of Theorem 2 is based on a recent renewal theorem for systems of equations reported in de Saporta (2003) and on an AR(1) recurrence equation satisfied by the discretization of $Y$ that we will define in Section 2. In Section 3 we study an operator related to our problem and prove Theorem 1. Sections 4–7 are devoted to the proof of Theorem 2. First we state two renewal theorems for systems of equations. Then in Section 5 we derive the renewal equations associated to our problem. In Sections 6 and 7 we prove Theorem 2, the latter section being dedicated to the proof that the constant $L$ is nonzero. Finally, in Section 8 we give further details on the computation of the exponent $\kappa$.

**2. Discretization of the process and an AR(1) equation.** First we give an explicit formula for the diffusion process. For $0 \le s \le t$, let

$$\Phi(s,t) = \Phi_{s,t}(\omega) = \exp\int_s^t a(X_u)\,du.$$

The process $Y$ has the representation [see Karatzas and Shreve (1991)]:

$$Y_t = Y_t(\omega) = \Phi(0,t)\left[Y_0 + \int_0^t \Phi(0,u)^{-1}\sigma(X_u)\,dW_u\right],$$

and for $0 \le s \le t$, $Y$ satisfies the recursion equation

$$Y_t = \Phi(s,t)\left[Y_s + \int_s^t \Phi(s,u)^{-1}\sigma(X_u)\,dW_u\right]$$

$$= \Phi(s,t)Y_s + \int_s^t\left[\exp\int_u^t a(X_v)\,dv\right]\sigma(X_u)\,dW_u.$$

It is useful to rewrite this recursion as

$$(3) \qquad Y_t(\omega) = \Phi_{s,t}(\omega)Y_s(\omega) + V_{s,t}^{1/2}(\omega)\xi_{s,t},$$

where $\xi_{s,t}$ is a standard Gaussian variable, function of $(W_u,\ s \le u \le t)$, and

$$V_{s,t}(\omega) = \int_s^t \exp\left[2\int_u^t a(X_v)\,dv\right]\sigma^2(X_u)\,du.$$

For $\delta > 0$, we will call *discretization at step size $\delta$ of $Y$* the discrete-time process $Y^{(\delta)} = (Y_{n\delta})_n$, where $n \in \mathbb{N}$. Our study of $Y$ is based on the investigations of these discretization $Y^{(\delta)}$ as in [Guyon, Iovleff and Yao](2004).

More precisely, for a fixed $\delta > 0$, the discretization $Y^{(\delta)}$ follows an AR(1) equation with random coefficients:

$$(4) \qquad Y_{(n+1)\delta}(\omega) = \Phi_{n+1}(\omega)Y_{n\delta}(\omega) + V_{n+1}^{1/2}(\omega)\xi_{n+1},$$

with

$$\Phi_{n+1}(\omega) = \Phi_{n+1}(\delta)(\omega) = \exp\left[\int_{n\delta}^{(n+1)\delta} a(X_u(\omega))\,du\right],$$

$$V_{n+1}(\omega) = \int_{n\delta}^{(n+1)\delta} \exp\left[2\int_u^{(n+1)\delta} a(X_v(\omega))\,dv\right]\sigma^2(X_u(\omega))\,du,$$

where $(\xi_n)$ is a standard Gaussian i.i.d. sequence defined on $(\Theta, \mathcal{B}, Q')$. Note that under Condition 2, all these discretizations are ergodic with the same limit distribution $\nu$ [see [Guyon, Iovleff and Yao](2004)].

## 3. Study of a related operator.
We now introduce a related operator $A$ and investigate its properties. Fix $s \geq 0$ and $\delta > 0$. We define the operator $A_{(s,\delta)}$ by

$$A_{(s,\delta)}\varphi(i) = \mathbb{E}_i[\Phi_1^s(\delta)\varphi(X_\delta)],$$

for every function $\varphi \colon E \to \mathbb{R}$ and every $i$ in $E$. It has the following semigroup property:

PROPOSITION 1. *Fix $s \geq 0$. Then for all $\delta, \gamma > 0$ we have*

$$A_{(s,\delta)}A_{(s,\gamma)} = A_{(s,\delta+\gamma)}.$$

PROOF. Set $\varphi \colon E \to \mathbb{R}$ and $i$ in $E$. We have

$$A_{(s,\delta)}A_{(s,\gamma)}\varphi(i) = \mathbb{E}_i[\Phi_1^s(\delta)A_{(s,\gamma)}\varphi(X_\delta)]$$
$$= \mathbb{E}_i[\Phi_1^s(\delta)\mathbb{E}_{X_\delta}[\Phi_1^s(\gamma)\varphi(X_\gamma)]]$$
$$= \mathbb{E}_i\left[\exp\left(s\int_0^\delta a(X_u)\,du\right)\mathbb{E}_{X_\delta}\left[\exp\left(s\int_0^\gamma a(X_u)\,du\right)\varphi(X_\gamma)\right]\right].$$

Then the Markov property yields

$$A_{(s,\delta)}A_{(s,\gamma)}\varphi(i) = \mathbb{E}_i\left[\exp\left(s\int_0^{\delta+\gamma} a(X_u)\,du\right)\varphi(X_{\delta+\gamma})\right]$$
$$= \mathbb{E}_i[\Phi_1^s(\delta+\gamma)\varphi(X_{\delta+\gamma})]$$
$$= A_{(s,\delta+\gamma)}\varphi(i). \qquad \square$$

Note that $A_{(s,\delta)}\varphi(i) = \sum_{j=1}^N \mathbb{E}_i[\Phi_1^s \mathbf{1}_{X_\delta = j}]\varphi(j)$, and therefore $A_{(s,\delta)}$ can be rewritten as the matrix $((A_{(s,\delta)})_{ij})_{1 \leq i,j \leq N}$ with $(A_{(s,\delta)})_{ij} = \mathbb{E}_i[\Phi_1^s \mathbf{1}_{X_\delta = j}]$. Note also that it is a positive operator.

3.1. *Spectral radius.* Now we investigate the properties of the spectral radius of $A$. First, we recall a result from Guyon, Iovleff and Yao (2004).

PROPOSITION 2. *Fix $s > 0$ and $\delta > 0$. Then $A_{(s,\delta)}$ is irreducible, aperiodic and satisfies*

$$(5) \qquad \mathbb{E}_\mu[(\Phi_1 \cdots \Phi_k)^s] = \sum_{i \in E} A_{(s,\delta)}^k \mathbf{1}(i)\mu(i) = \mu A_{(s,\delta)}^k \mathbf{1},$$

*where $\mathbf{1}$ is the constant function equal to $1$ on $E$.*

We denote by $\rho(X)$ the spectral radius of a matrix $X$. Proposition 2 yields the following corollaries.

COROLLARY 1. *We have*

$$\rho(A_{(s,\delta)}) = \lim_{k \to \infty} (\mathbb{E}_\mu[(\Phi_1 \cdots \Phi_k)^s])^{1/k}.$$

PROOF. As $A_{(s,\delta)}$ is a (component-wise) positive, irreducible and aperiodic matrix, Theorem 8.5.1 of Horn and Johnson (1985) gives the existence of a matrix $B_{(s,\delta)}$ with positive coefficients such that

$$(6) \qquad \frac{(A_{(s,\delta)})^n}{(\rho(A_{(s,\delta)}))^n} \underset{n \to \infty}{\longrightarrow} B_{(s,\delta)}.$$

This result and (5) yield the expected result. □

COROLLARY 2. *For all fixed $\delta > 0$, the mapping $s \longmapsto \log \rho(A_{(s,\delta)})$ is convex on $\mathbb{R}_+$.*

Note that for all fixed $\delta > 0$ and $i$ in $E$, we have $A_{(0,\delta)}\mathbf{1}(i) = \mathbb{E}_i(1) = 1$. Thus, as $A_{(0,\delta)}$ is a positive operator, it is also a stochastic matrix and $\rho(A_{(0,\delta)}) = 1$.

PROPOSITION 3. *For all fixed $\delta > 0$, the right-hand derivative of the mapping $s \longmapsto \log \rho(A_{(s,\delta)})$ at $0$ is negative.*

PROOF. As all the functions considered are convex, we have

$$\frac{\partial}{\partial s} \log(\rho(A_{(s,\delta)})) = \lim_{n \to \infty} \frac{\partial}{\partial s} \frac{1}{n} \log \mathbb{E}_\mu[(\Phi_1 \cdots \Phi_n)^\kappa]$$

$$= \lim_{n \to \infty} \frac{1}{n} \frac{\mathbb{E}_\mu[(\Phi_1 \cdots \Phi_n)^\kappa \cdot \sum_{i=1}^n \log \Phi_i]}{\mathbb{E}_\mu[(\Phi_1 \cdots \Phi_n)^\kappa]}.$$

The sequence $(\Phi_n)$ is stationary, thus the ergodic theorem yields

(7) $$\frac{1}{n}\sum_{k=1}^{n}\log\Phi_k \underset{n\to\infty}{\longrightarrow} \mathbb{E}_\mu[\log\Phi_1], \qquad \mathbb{P}_\mu\text{-almost surely.}$$

But $\mathbb{E}_\mu[\log\Phi_1] < 0$ because of Condition 2. Indeed, we have

$$\mathbb{E}_\mu[\log\Phi_1] = \mathbb{E}_\mu\left[\int_0^\delta a(X_u)\,du\right] = \int_0^\delta \mathbb{E}_\mu[a(X_u)]\,du = \delta\alpha < 0.$$

Thus we get, as expected,

$$\frac{\partial}{\partial s}\bigg|_{s=0}\log\left(\rho(A_{(s,\delta)})\right) = \lim_{n\to\infty}\frac{1}{n}\mathbb{E}_\mu\left[\sum_{i=1}^{n}\log\Phi_i\right]$$
$$= \mathbb{E}_\mu[\log\Phi_1] < 0. \qquad \square$$

COROLLARY 3.   *Fix $\delta > 0$. We have the following dichotomy:*

(i)  *either for all $s > 0$, $\rho(A_{(s,\delta)}) < 1$,*
(ii)  *or there exists a unique $\kappa > 0$ such that $\rho(A_{(\kappa,\delta)}) = 1$, and in this case $\rho(A_{(s,\delta)}) > 1$ for all $s > \kappa$ and $\rho(A_{(s,\delta)}) < 1$ for all $0 < s < \kappa$.*

3.2. *Choice of $\delta$.*  Now we are going to prove that the preceding dichotomy is in fact independent of the value of $\delta$.

PROPOSITION 4.   *Fix $s \geq 0$. The following propositions are equivalent:*

(i)  *there exists $\delta > 0$ such that $\rho(A_{(s,\delta)}) < 1$,*
(ii)  *for all $\delta > 0$ we have $\rho(A_{(s,\delta)}) < 1$.*

*The same equivalence is true if we replace "$< 1$" by "$> 1$" or "$= 1$."*

PROOF.   Set $\delta > 0$ such that $\rho(A_{(s,\delta)}) < 1$, and $\gamma > 0$. For all integer $n \geq 1$ we define $m_n \in \mathbb{N}^*$ and $0 \leq \beta_n < \delta$ by $n\gamma = m_n\delta + \beta_n$ ($m_n$ the integer part of $n\gamma/\delta$ and $\beta_n$ its fractional part multiplied by $\delta$). Thus Proposition 1 yields

$$A^n_{(s,\gamma)} = A_{(s,n\gamma)} = A^{m_n}_{(s,\delta)}A_{(s,\beta_n)}.$$

But for all $n$ we have

$$\|A_{(s,\beta_n)}\| \leq \max_i \mathbb{E}_i[\Phi_1^s(\beta_n)]$$
$$\leq \exp\left(s\beta_n\max_i(a_i)\right)$$
$$\leq \exp\left(s\delta\max_i(a_i)\right).$$

This upper bound is independent of $n$. Thus we have

$$\log \|A_{(s,\gamma)}^n\| \leq \log \|A_{(s,\delta)}^{m_n}\| + c,$$

where $c$ is a positive constant. We get

$$\log \rho(A_{(s,\gamma)}) = \lim_n \frac{1}{n} \log \|A_{(s,\gamma)}^n\|$$

$$\leq \limsup_n \frac{1}{n} \log \|A_{(s,\delta)}^{m_n}\|$$

$$= \frac{\gamma}{\delta} \log \rho(A_{(s,\delta)}),$$

as $m_n \sim n\gamma\delta^{-1}$. Hence $\rho(A_{(s,\gamma)}) \leq \rho(A_{(s,\delta)})^{\gamma/\delta} < 1$.

For the case "$=1$," fix $\delta_0$ and a corresponding $\kappa$ such that $\rho(A_{(\kappa,\delta_0)}) = 1$. The mapping $s \longmapsto \rho(A_{(s,\delta_0)})$ is log-convex hence continuous. Thus we have

$$\rho(A_{(\kappa,\delta_0)}) = \sup_{s < \kappa} \rho(A_{(s,\delta_0)}).$$

Set $\delta > 0$. We want to prove that $\rho(A_{(\kappa,\delta)}) = 1$. According to Corollary 3, for all $s < \kappa$ we have $\rho(A_{(s,\delta_0)}) < 1$. Thus the preceding study yields that for all $s < \kappa$ we also have $\rho(A_{(s,\delta)}) < 1$. Hence we have

$$\rho(A_{(\kappa,\delta)}) = \sup_{s < \kappa} \rho(A_{(s,\delta)}) \leq 1.$$

Suppose that $\rho(A_{(\kappa,\delta)}) < 1$; then the first case implies again that $\rho(A_{(\kappa,\delta_0)}) < 1$, which is impossible. Thus we have $\rho(A_{(\kappa,\delta)}) = 1$ as expected.

The case "$>1$" is a consequence of these two cases and Corollary 3.  □

In the following we will write $A_s$ instead of $A_{(s,\delta)}$ each time it is nonambiguous. We have an easy criterion to know in which case we are.

PROPOSITION 5.  *The following properties are equivalent:*

  (i) *for all $i$ in $E$, $a(i) \leq 0$,*
  (ii) *for all $s > 0$, $\rho(A_s) < 1$.*

PROOF.  Suppose that for all $i$ in $E$ we have $a(i) \leq 0$. Fix $\delta > 0$. Then for all $s > 0$, we have $\Phi_1^s \leq 1$. Thus for all $i$, $A_s \mathbf{1}(i) = \mathbb{E}_i[\Phi_1^s] \leq 1$, and componentwise we have $A_s \mathbf{1} \leq \mathbf{1}$, which implies that $\rho(A_s) \leq 1$ for all $s > 0$. Corollary 3 then yields that for all $s$, we have actually $\rho(A_s) < 1$.

Now suppose there exists an $i_0$ such that $a(i_0) > 0$. Fix $s \geq 2\lambda(i_0)a(i_0)^{-1}$. It is proved in Guyon, Iovleff and Yao (2004) that for all function $\varphi$ from $E$ into $\mathbb{R}$ and all $i$ in $E$ we have for small $\delta$,

$$(8) \qquad A_s\varphi(i) = [1 + \delta(sa(i) - \lambda(i))]\varphi(i) + \delta\lambda(i) \sum_{j \neq i}[q(i,j)\varphi(j)] + o(\delta).$$

Let $\psi$ be the function from $E$ into $\mathbb{R}$ such that $\psi(i_0) = 1$ and $\psi(i) = 0$ for all $i \neq i_0$. Then for all $i \neq i_0$ we have $A_s \psi(i) = \mathbb{E}_i[\Phi_1^s \mathbf{1}_{X_\delta = i_0}] \geq 0$ and for $i = i_0$ we have

$$A_s \psi(i_0) = 1 + \delta(sa(i_0) - \lambda(i)) + o(\delta) \geq 1 + \delta \frac{sa(i_0)}{2} + o(\delta)$$

as we have chosen $s \geq 2\lambda(i)a(i_0)^{-1}$. Thus component-wise, for small enough $\delta$, we have

$$A_s \psi \geq \left( 1 + \delta \frac{sa(i_0)}{2} + o(\delta) \right) \psi$$

$$\geq \left( 1 + \delta \frac{sa(i_0)}{4} \right) \psi.$$

Thus $\rho(A_s) \geq 1 + \delta \frac{sa(i_0)}{4} > 1.$   $\square$

This proposition ends the proof of Theorem 1 since we have the following result from Guyon, Iovleff and Yao (2004) that relates the spectral radius of $A_s$ to the moments of the stationary law $\nu$:

PROPOSITION 6.  *Set* $s > 0$. *If* $\rho(A_s) < 1$, *then the stationary law* $\nu$ *of* $Y$ *has a moment of order* $s$.

The proof of Theorem 1 is now complete.

**4. Renewal theory for systems.**  Now we proceed to prove Theorem 2. From now on, we will assume that there is an $i$ such that $a(i) > 0$. Our approach is based on a new renewal theorem for systems of renewal equations. First we introduce some notation and conventions that we will apply throughout.

Let $F = (F_{ij})_{1 \leq i,j \leq p}$ be a matrix of distributions: nondecreasing, right-continuous functions on $\mathbb{R}$ into $\mathbb{R}_+$ with limit 0 at $-\infty$. For all $p \times r$ matrix $H$ of Borel-measurable, real-valued functions $H_{ij}$ on $\mathbb{R}$ that are bounded on compact intervals, we define the *convolution product* $F * H$ by

$$(F * H)_{ij}(t) = \sum_{k=1}^p \int_{-\infty}^\infty H_{kj}(t - u) F_{ik}(du)$$

where it exists.

The transpose of a vector or matrix $X$ will always be denoted $X'$. We study the renewal equation $Z = F * Z + G$, where $G = (G_1, \ldots, G_p)'$ is a vector of Borel-measurable, real-valued functions, bounded on compact intervals, and $Z = (Z_1, \ldots, Z_p)'$ is a vector of functions. The renewal theorem will give the limit of $Z$ at $+\infty$.

For all real $t$, we set:

(a) $B = (b_{ij})_{1 \le i,j \le p}$ where $b_{ij} = \int u F_{ij}(du)$ if it exists, the expectation of $F$,

(b) $F^{(0)}(t) = (\delta_{ij}(t))_{1 \le i,j \le p}$ where $\delta_{ij}(t) = \mathbf{1}_{t \ge 0}$ if $i = j$ and $0$ otherwise, so that $F^{(0)} * H = H$ for all $H$ as in the definition above,

(c) $F^{(n)}(t) = F * F^{(n-1)}(t)$, the $n$-fold convolution of $F$,

(d) $U(t) = \sum_{n=0}^{\infty} F^{(n)}(t)$, the *renewal function* associated with $F$.

We will also assume that all the measures $F_{ij}$ are finite:

$$F_{ij}(\infty) = \lim_{t \to \infty} F_{ij}(t) < \infty,$$

and that $F(\infty)$ is an irreducible matrix. $F(\infty)$ being an irreducible nonnegative matrix, we can use the Perron–Frobenius theorem: its spectral radius $\rho(F(\infty))$ is a simple eigenvalue with right and left positive eigenvectors. We will also assume that $\rho(F(\infty)) = 1$, and we choose two positive eigenvectors $m$ and $u$ so that

$$F(\infty)m = m, \qquad u'F(\infty) = u', \qquad \sum_{i=1}^{p} m_i = 1, \qquad \sum_{i=1}^{p} u_i m_i = 1.$$

We also assume that the sequence $(\|F(\infty)^n\|)$ is bounded [e.g., if $F(\infty)$ is aperiodic, this is true]. We recall the following definition: $F$ is *lattice* if the following conditions are satisfied:

(a) For all $i \ne j$, $F_{ij}$ is concentrated on a set of the form $b_{ij} + \lambda_{ij}\mathbb{Z}$.

(b) For all $i$, $F_{ii}$ is concentrated on a set of the form $\lambda_{ii}\mathbb{Z}$.

(c) Each $\lambda_{ii}$ is an integral multiple of the same number. We take $\lambda$ to be the largest such number.

(d) For all $a_{ij}$, $a_{jk}$, $a_{ik}$ points of increase of $F_{ij}$, $F_{jk}$, $F_{ik}$, respectively, $a_{ij} + a_{jk} - a_{ik}$ is an integral multiple of $\lambda$.

We can now state the renewal theorem. It extends a previous result of Crump (1970) and Athreya and Rama Murthy (1976) which deals with the positive case: each distribution $F_{ij}$ has support on $\mathbb{R}_+$. The proof of this theorem is given in de Saporta (2003).

RENEWAL THEOREM A. *Assume that $F$ is as above and that, in addition, it is a nonlattice matrix, that its expectation $B$ exists, and that for all $t \in \mathbb{R}$, $U(t)$ is finite. If $G$ is directly Riemann integrable [see Feller (1966)], and $Z = U * G$ exists, then for all $i$, we have*

$$\lim_{t \to \infty} Z_i(t) = c m_i \sum_{j=1}^{p} \left[ u_j \int_{-\infty}^{\infty} G_j(y) \, dy \right],$$

*where $m$ and $u$ are the eigenvectors defined above and $c = (u'Bm)^{-1}$ (under these assumptions, $u'Bm \ne 0$).*

We also recall Theorem 2.3 of Athreya and Rama Murthy (1976) that will be used in Section 7. Note that this theorem can now be seen as a corollary of Theorem A.

RENEWAL THEOREM B.  *Let $F$ be a nonlattice matrix of distributions with support on the positive half-line, such that:*

(i) $\rho(F(0)) < 1$,
(ii) $F(\infty)$ *is finite, irreducible and aperiodic,*
(iii) *there exist $i$ and $j$ such that $F_{ij}(0) < F_{ij}(\infty)$.*

*Assume also that there is an $\alpha > 0$ such that $\rho(F_\alpha) = 1$, where $(F_\alpha)_{ij} = \int_0^\infty e^{-\alpha u} F_{ij}(du)$. Then for all $h > 0$, and all $i, j$, we have*

$$\lim_{t \to \infty} \int_t^{t+h} e^{-\alpha y} U_{ij}(dy) = c m_i u_j h,$$

*where $m$ and $u$ are right and left eigenvectors of $F_\alpha$, with the same normalization as above, $c = (u'Bm)^{-1}$, and $B = b_{ij}$ with $b_{ij} = \int_0^\infty u e^{-\alpha u} F_{ij}(du)$, $c$ being interpreted as zero if some $b_{ij}$ is equal to infinity.*

**5. The renewal equations.**  Now we are going to derive the renewal equations associated to our problem. In the following, we will suppose that the assumptions of Theorem 2 are satisfied. We set $\delta = 1$, and $\kappa$ will denote the unique positive solution of $\rho(A_s) = 1$. We are going to study the discretization $Y^{(1)}$.

5.1. *Notation.*  As $X$ is a stationary process, we can extend it to negative $t$ and define the coefficients $\Phi_n$, $V_n$ and $\xi_n$ for negative values of $n$. Let $b_n = V_n^{1/2} \xi_n$ and

$$R_n = \sum_{k=0}^\infty \Phi_n \Phi_{n-1} \cdots \Phi_{n-k+1} b_{n-k}$$

(instead of $\widetilde{Y}_n$) be the unique stationary solution of (4): $R_{n+1} = \Phi_{n+1} R_n + b_{n+1}$. The limit law $\nu$ of $Y$ is also the law of $R_1$. Thus we are going to study the random variable $R_1$.

The tail of the stationary solution of such recursive equations has already been studied in various cases. In the i.i.d. multidimensional case: $\Phi_n$ are matrices and $R_n$ and $b_n$ vectors, renewal theory is used in Kesten (1973) to prove a heavy tail property when the $\Phi_n$ either have a density or are nonnegative. These results were extended in Le Page (1983) to a wider class of i.i.d. random matrices. Finally, in Goldie (1991) a new specific implicit renewal theorem is proved and the same results are derived in the i.i.d. one-dimensional case. This theorem also applies to the study of the tail of several other random recurrences implying i.i.d. random variables. Recently, Goldie's results were extended in de Saporta (2004) to the case where $(\Phi_n)$ is a finite state space Markov chain. Here, $(\Phi_n)$ is not a Markov chain, but conditionally to $X_n$, $\Phi_n$ and $\Phi_{n+1}$ are independent. Our proof is thus very similar to that of de Saporta (2004), but we will repeat all the details for completeness.

Note that $\xi_n$ are standard Gaussian random variables, thus they are symmetric, and they are also independent from the sequences $(\Phi_n)$ and $(V_n)$. Hence we have

$$\mathbb{P}_\mu\left(\sum_{k=0}^\infty \Phi_1\Phi_0\cdots\Phi_{2-k}b_{1-k} > t\right)$$

$$= \mathbb{P}_\mu\left(\sum_{k=0}^\infty \Phi_1\Phi_0\cdots\Phi_{2-k}V_{1-k}^{1/2}\xi_{1-k} > t\right)$$

$$= \mathbb{P}_\mu\left(\sum_{k=0}^\infty \Phi_1\Phi_0\cdots\Phi_{2-k}V_{1-k}^{1/2}(-\xi_{1-k}) > t\right)$$

$$= \mathbb{P}_\mu\left(-\sum_{k=0}^\infty \Phi_1\Phi_0\cdots\Phi_{2-k}b_{1-k} > t\right).$$

Thus we have $\nu(]t, +\infty[) = \nu(]-\infty, -t[)$ for all $t$; hence if one of the limits stated in Theorem 2 exists, the other exists too and equals the same value. Therefore we need study only one limit.

To study the tail of $R_1$, we introduce a new function. For all $t$ in $\mathbb{R}$, we set

$$z(t) = e^{-t}\int_0^{e^t} u^\kappa \mathbb{P}(R_1 > u)\, du.$$

Lemma 9.3 of Goldie (1991) ascertains that if $z(t)$ has a limit when $t$ tends to infinity, then $t^\kappa \mathbb{P}(R_1 > t)$ also has the same limit.

For all $i$ in $E$ and $t$ in $\mathbb{R}$, we also set

$$Z_i(t) = e^{-t}\int_0^{e^t} u^\kappa \mathbb{P}(R_1 > u, X_1 = i)\, du,$$

so that $z(t) = \sum_{i=1}^N Z_i(t)$. We are now going to prove that $Z = {}^t(Z_1, \ldots, Z_N)$ satisfies a system of renewal equations.

5.2. *The renewal equations.* As $R_n$ satisfies (4), we have $R_1 = \Phi_1 R_0 + b_1$; thus for all $t$ in $\mathbb{R}$, we have

$$\mathbb{P}_\mu(R_1 > u, X_1 = i) = \mathbb{P}_\mu(\Phi_1 R_0 > u, X_1 = i) + \psi_i(u),$$

where

$$\psi_i(t) = \mathbb{P}_\mu(t - b_1 < \Phi_1 R_0 \leq t, X_1 = i) - \mathbb{P}_\mu(t < \Phi_1 R_0 \leq t - b_1, X_1 = i).$$

We set $G_i(t) = e^{-t} \int_0^{e^t} u^\kappa \psi_i(u)\, du$, and $G = {}^t(G_1, \ldots, G_N)$. Then we have

$$z(t) = \sum_{i=1}^N \left[ e^{-t} \int_0^{e^t} u^\kappa \mathbb{P}_\mu(\Phi_1 R_0 > u, X_1 = i)\, du + G_i(t) \right].$$

We have $\Phi_1 \geq 0$ and conditionally to $X_0$, $\Phi_1$ and $R_0$ are independent. Thus, a simple change of variable and stationarity yield

$$e^{-t} \int_0^{e^t} u^\kappa \mathbb{P}_\mu(\Phi_1 R_0 > u, X_1 = i)\, du$$

$$= \sum_{j=1}^N e^{-t} \int_0^{e^t} u^\kappa \mathbb{P}_\mu(\Phi_1 R_0 > u, X_1 = i | X_0 = j)\mu(j)\, du$$

$$= \sum_{j=1}^N e^{-t} \int_0^{e^t} u^\kappa \mathbb{P}_j(\Phi_1 R_0 > u, X_1 = i)\mu(j)\, du$$

$$= \sum_{j=1}^N \mathbb{E}_j \left[ \Phi_1^\kappa \mathbf{1}_{X_1=i} e^{-(t-\log \Phi_1)} \int_0^{e^{t-\log \Phi_1}} u^\kappa \mathbb{P}_j(R_0 > u)\, du \right] \mu(j)$$

$$= \sum_{j=1}^N \mathbb{E}_j \left[ \Phi_1^\kappa \mathbf{1}_{X_1=i} e^{-(t-\log \Phi_1)} \int_0^{e^{t-\log \Phi_1}} u^\kappa \mathbb{P}_\mu(R_0 > u | X_0 = j)\, du \right] \mu(j)$$

$$= \sum_{j=1}^N \mathbb{E}_j \left[ \Phi_1^\kappa \mathbf{1}_{X_1=i} e^{-(t-\log \Phi_1)} \int_0^{e^{t-\log \Phi_1}} u^\kappa \mathbb{P}_\mu(R_1 > u, X_1 = j)\, du \right].$$

Thus we get the following system of equations: for all $i$ in $E$, we have

$$Z_i(t) = \sum_{j=1}^N \left[ \mathbb{E}_j[\Phi_1^\kappa \mathbf{1}_{X_1=i} Z_j(t - \log \Phi_1)] \right] + G_i(t)$$

(9)

$$= \sum_{j=1}^N [F_{ij} * Z_j(t)] + G_i(t),$$

where $F_{ij}(t) = \mathbb{E}_j[\Phi_1^\kappa \mathbf{1}_{X_1=i} \mathbf{1}_{t \geq \log \Phi_1}]$. Thus $F = (F_{ij})_{i,j \in E}$ is a matrix of distributions in the sense of Section 4, and system (9) is a system of renewal equations that can be rewritten as $Z = F * Z + G$. To apply Theorem A, we now have to prove that $F$ and $G$ satisfy its assumptions.

**6. Proof of Theorem 2, part I.** As $E$ is a finite set, $\Phi_1$ is bounded. Therefore, for all $i,j$ in $E$, the measures $F_{ij}$ are finite and $F_{ij}(\infty) = \mathbb{E}_j[\Phi_1^\kappa \mathbf{1}_{X_1 = i}]$. Note that $F(\infty) = A'_\kappa$. As $A_\kappa$ is irreducible and aperiodic by Proposition 2, so is $F(\infty)$, and its spectral radius also equals to 1. Besides, we have $b_{ij} = \mathbb{E}_j[\Phi_1^\kappa \mathbf{1}_{X_1 = i} \log \Phi_1]$, thus the $F_{ij}$ have finite expectation.

We are going to prove that the other assumptions of Theorem A are valid here in the following sections.

6.1. *$F$ is nonlattice.* Set $a_m = \min_{i \in E}\{a(i)\}$, $a_M = \max_{i \in E}\{a(i)\}$ and $i_0$, $j_0$ in $E$ such that $a(i_0) = a_m$ and $a(j_0) = a_M$.

PROPOSITION 7. *For all $i,j$ in $E$, $x \in ]a_m, a_M[$ and small enough $\varepsilon > 0$, we have*

$$\mathbb{P}_i\left( \int_0^1 a(X_u)\,du \in ]x - \varepsilon; x + \varepsilon[, \ X_1 = j \right) > 0,$$

*that is, $x$ is a point of increase of $\log \Phi_1$ conditionally to $X_0 = i$ and $X_1 = j$.*

PROOF. Set $x \in ]a_m, a_M[$ and $0 < t < 1$ such that $x = ta_m + (1-t)a_M$. Fix $i$ and $j$ in $E$. As $q$ is an irreducible matrix, we can find integers $0 \le l \le m \le n$ and $k_1, \ldots, k_n$ in $E$ such that $q_{i,k_1}q_{k_1,k_2}\cdots q_{k_l,i_0} > 0$, $q_{i_0,k_{l+1}}q_{k_{l+1}k_{l+2}}\cdots q_{k_m,j_0} > 0$ and $q_{j_0,k_{m+1}}q_{k_{m+1}k_{m+2}}\cdots q_{k_n,j} > 0$. Set also $y = a(i) + a(k_1) + \cdots + a(k_l) - (l+1)a_m + a(k_{l+1}) + \cdots + a(k_m) - (n-l+1)a_M + a(k_{m+1}) + \cdots + a(k_n) + a(j)$, and $z = \min\{\varepsilon|y|^{-1}, t(l+1)^{-1}, (1-t)(n-l+1)^{-1}\}$. Then we have

$$\mathbb{P}_i\left( \int_0^1 a(X_u)\,du \in ]x - \varepsilon; x + \varepsilon[, \ X_1 = j \right)$$

$$\ge \mathbb{P}_i(X_u = i \text{ on } [0; \eta[, \ X_u = k_1 \text{ on } [\eta; 2\eta[, \ldots, \ X_u = k_l \text{ on } [l\eta; (l+1)\eta[,$$

$$X_u = i_0 \text{ on } [(l+1)\eta, t[, \ X_u = k_{l+1} \text{ on } [t; t+\eta[, X_u = k_{l+2} \text{ on }$$

(10) $$[t+\eta; t+2\eta[, \ldots, \ X_u = k_m \text{ on } [t+(m-l-1)\eta; t+(m-l)\eta[,$$

$$X_u = j_0 \text{ on } [t+(m-l)\eta; 1-(n-m+1)\eta[, \ X_u = k_{m+1} \text{ on }$$

$$[1-(n-m+1)\eta; 1-(n-m)\eta[, \ldots, \ X_u = k_n \text{ on } [1-2\eta; 1-\eta[,$$

$$X_u = j \text{ on } [1-\eta; 1]; \ \eta \in ]0; z[).$$

Indeed, on this event we have

$$\int_0^1 a(X_u)\,du$$

$$= \eta a(i) + \eta a(k_1) + \cdots + \eta a(k_l) + (t - (l+1)\eta)a_m + \eta a(k_{l+1})$$

$$+ \cdots + \eta a(k_m) + ((1-t) - (n-l+1)\eta)a_M + \eta a(k_{m+1})$$

$$+ \cdots + \eta a(k_n) + \eta a(j)$$

$$= ta_m + (1-t)a_M + \eta y$$

$$= x + \eta y,$$

thus if $\eta < \varepsilon |y|^{-1}$, then we have $\int_0^1 a(X_u)\,du \in\, ]x-\varepsilon; x+\varepsilon[$. Probability (10) can be computed [see, e.g., Norris (1998)]:

$$(10) = \mu(i)q_{i,k_1}q_{k_1,k_2}\cdots q_{k_l,i_0}q_{i_0,k_{l+1}}\cdots q_{k_m,j_0}q_{j_0,k_{m+1}}\cdots q_{k_n,j}$$

$$\times \lambda(i)\lambda(k_1)\cdots\lambda(k_n)\lambda(i_0)(l-1)\lambda(j_0)(n-l+1)$$

$$\times \int_0^z [e^{-\lambda(i)\eta}e^{-\lambda(k_1)\eta}$$

$$\cdots e^{-\lambda(k_n)\eta}e^{-\lambda(i_0)(t-(l-1))\eta}e^{-\lambda(j_0)(1-t-(n-l+1)\eta}e^{-\lambda(j)\eta}]\,d\eta.$$

Thus our choice of $k_1, \ldots, k_n$ and $z$ ascertains that this probability is positive, which proves the proposition.  □

Therefore none of the $F_{ij}(\cdot) = \mathbb{E}_j[\Phi_1^\kappa \mathbf{1}_{X_1=i}\mathbf{1}_{\cdot \geq \log \Phi_1}]$ can be concentrated on a lattice set, and in particular $F$ is nonlattice.

6.2. *Finiteness of $U$.*  We are going to prove that for all $i, j$ in $E$ and $t$ in $\mathbb{R}$, $U_{ij}(t)$ is finite. We start with computing the $n$-fold convolution of $F$.

LEMMA 1.  *For all $n, i, j, t$ we have*

$$F_{ij}^{(n)}(t) = \mathbb{E}_j[\Phi_1^\kappa \cdots \Phi_n^\kappa \mathbf{1}_{\log \Phi_1 \cdots \Phi_n \geq t}\mathbf{1}_{X_n=i}].$$

PROOF.  For $n=1$, it is the definition of $F$. Suppose the formula is true for a fixed $n$. Then the Markov property and stationarity yield

$$F_{ij}^{(n+1)}(t)$$

$$= \sum_{k=1}^N F_{ik} * F_{kj}^{(n)}(t) = \sum_{k=1}^N \int F_{kj}^{(n)}(t-u)F_{ik}(du)$$

$$= \sum_{k=1}^N \int \mathbb{E}_j[\Phi_1^\kappa \cdots \Phi_n^\kappa \mathbf{1}_{\log \Phi_1 \cdots \Phi_n \leq t-u}\mathbf{1}_{X_n=k}]\mathbb{E}_k[\Phi_1^\kappa \delta_u(\log \Phi_1)\mathbf{1}_{X_1=i}]$$

$$= \sum_{k=1}^N \int \mathbb{E}_\mu[\Phi_1^\kappa \cdots \Phi_n^\kappa \mathbf{1}_{\log \Phi_1 \cdots \Phi_n \leq t-u}\mathbf{1}_{X_n=k}\mathbf{1}_{X_0=j}]$$

$$\times \mathbb{E}_\mu[\Phi_{n+1}^\kappa \delta_u(\log \Phi_{n+1}) \mathbf{1}_{X_{n+1}=i} \mathbf{1}_{X_n=k}] \frac{1}{\mu(k)\mu(j)}$$

$$= \sum_{k=1}^N \mathbb{E}_\mu[\Phi_1^\kappa \cdots \Phi_n^\kappa \mathbf{1}_{\log \Phi_1 \cdots \Phi_n \le t - \log \Phi_{n+1}} \mathbf{1}_{X_0=j} | \mathbf{1}_{X_n=k}]$$

$$\times \mathbb{E}_\mu[\Phi_{n+1}^\kappa \mathbf{1}_{X_{n+1}=i} | \mathbf{1}_{X_n=k}] \frac{\mu(k)}{\mu(j)}$$

$$= \sum_{k=1}^N \mathbb{E}_\mu[\Phi_1^\kappa \cdots \Phi_n^\kappa \Phi_{n+1}^\kappa \mathbf{1}_{\log \Phi_1 \cdots \Phi_n \le t - \log \Phi_{n+1}} \mathbf{1}_{X_0=j} \mathbf{1}_{X_{n+1}=i} | \mathbf{1}_{X_n=k}] \frac{\mu(k)}{\mu(j)}$$

$$= \mathbb{E}_\mu[\Phi_1^\kappa \cdots \Phi_n^\kappa \Phi_{n+1}^\kappa \mathbf{1}_{\log \Phi_1 \cdots \Phi_n \Phi_{n+1} \le t} \mathbf{1}_{X_0=j} \mathbf{1}_{X_{n+1}=i}] \frac{1}{\mu(j)}$$

$$= \mathbb{E}_j[\Phi_1^\kappa \cdots \Phi_n^\kappa \Phi_{n+1}^\kappa \mathbf{1}_{\log \Phi_1 \cdots \Phi_n \Phi_{n+1} \le t} \mathbf{1}_{X_{n+1}=i}].$$

Thus the formula is also true for $n+1$ and the lemma is proved. $\square$

We have seen that $F(\infty) = A_\kappa'$. Proposition 1 and the preceding lemma also imply that for all $n$ we have $F^{(n)}(\infty) = (A_\kappa^n)' = F(\infty)^n$. We can prove a more general result.

Lemma 2. *For all $n$ and $0 \le r < \kappa$ we have*

$$\int_{-\infty}^\infty e^{-ru} F^{(n)}(du) = (A_{\kappa-r}^n)'.$$

Proof. For all $i, j$ in $E$, Proposition 1 and the preceding lemma yield

$$\int_{-\infty}^\infty e^{-ru} F_{ij}^{(n)}(du)$$

$$= \int_{-\infty}^\infty e^{-ru} \mathbb{E}_j[\Phi_1^\kappa \cdots \Phi_n^\kappa \delta_u(\log \Phi_1 \cdots \Phi_n) \mathbf{1}_{X_n=i}]$$

$$= \mathbb{E}_j[\Phi_1^\kappa \cdots \Phi_n^\kappa e^{-r \log \Phi_1 \cdots \Phi_n} \mathbf{1}_{X_n=i}]$$

$$= \mathbb{E}_j[\Phi_1^{\kappa-r} \cdots \Phi_n^{\kappa-r} \mathbf{1}_{X_n=i}]$$

$$= (A_{\kappa-r}^n)_{ji}. \qquad \square$$

Now fix $0 < r < \kappa$. We have

$$
\begin{aligned}
U_{ij}(t) = \sum_{n=0}^\infty F_{ij}^{(n)}(t) &\le e^{rt} \int_{-\infty}^t e^{-ru} \sum_{n=0}^\infty F_{ij}^{(n)}(du) \\
&\le e^{rt} \sum_{n=0}^\infty \int_{-\infty}^\infty e^{-ru} F_{ij}^{(n)}(du) = e^{rt} \sum_{n=0}^\infty (A_{\kappa-r}^n)_{ji},
\end{aligned}
$$
(11)

according to the preceding lemma. But Corollary 3 says that $\rho(A_{\kappa-r}) < 1$. Thus the series in (11) converges. Hence $U_{ij}(t) < \infty$ for all $i, j$ in $E$ and $t$ in $\mathbb{R}$.

6.3. *Proof of* $Z = U * G$. Iterating the renewal equation (9) yields, for all $n$,

$$(12) \qquad Z = F^{(n)} * Z + \sum_{k=0}^{n-1} F^{(k)} * G.$$

The same change of variable as in Section 5.2 yields

$$\sum_{i=1}^{N} (F^{(n)} * Z)_i(t) = e^{-t} \int_0^{e^t} u^\kappa \mathbb{P}_\mu(\Phi_1 \Phi_2 \cdots \Phi_n R_0 > u) \, du.$$

But we have seen at (7) that we have $\Phi_1 \cdots \Phi_n \to 0$ when $n$ tends to infinity. Thus the bounded convergence Theorem yields $\sum_{i=1}^{N}(F^{(n)} * Z)(x,t) \to 0$ as $n$ tends to infinity. Each term of this sum is nonnegative, thus each term tends to 0. Letting $n$ tend to infinity in (12), we thus get $Z = U * G$.

6.4. *$G$ is directly Riemann integrable.* As the $G_i$ are clearly continuous in $t$, it is sufficient to prove that

$$\sum_{l=-\infty}^{\infty} \sup_{l \le t < l+1} |G_i(t)| < \infty$$

[see Feller (1966)]. But for all $i, t$, we have $G_i(t) = G_i^1(t) - G_i^2(t)$, where

$$G_i^1(t) = e^{-t} \int_0^{e^t} u^\kappa \mathbb{P}_\mu(u - b_1 < \Phi_1 R_0 \le u, X_1 = i) \, du \ge 0,$$

$$G_i^2(t) = e^{-t} \int_0^{e^t} u^\kappa \mathbb{P}_\mu(u < \Phi_1 R_0 \le u - b_1, X_1 = i) \, du \ge 0.$$

For all real $t$, we have $G_i(t) \le G_i^1(t) \le e^{-t} \int_0^{e^t} u^\kappa \, du = e^{t(\kappa+1)}(\kappa+1)^{-1}$. In particular, $G_i$ is directly Riemann integrable on $\mathbb{R}_-$. We still have to study $G_i^1$ and $G_i^2$ on $\mathbb{R}_+$. These two functions being of the same kind, we only give the detailed study of the first one.

The proof is adapted from Le Page (1983). Set $\varepsilon \in {]0; 1[}$ such that $-1 < \kappa - (1-\varepsilon) < 0$. Thus we have

$$
\begin{aligned}
(13) \qquad 0 \le e^t G_i^1(t) &\le \int_0^{e^t} u^\kappa \mathbb{P}_\mu(b_1 > u^\varepsilon, X_1 = i) \, du \\
&\quad + \int_0^{e^t} u^\kappa \mathbb{P}_\mu(u - u^\varepsilon < \Phi_1 R_0 \le u, X_1 = i) \, du.
\end{aligned}
$$

We are going to give an upper bound for each one of these two terms.

*First term.* Chebychev inequality yields

$$(14) \qquad \int_0^{e^t} u^\kappa \mathbb{P}_\mu(b_1 > u^\varepsilon, X_1 = i)\, du \le \mathbb{E}_\mu |b_1|^\kappa \frac{e^{t(1+\kappa(1-\varepsilon))}}{1 + \kappa(1-\varepsilon)}.$$

Note that $b_1$ has moments of all order. Indeed, we have, by independence, $\mathbb{E}_\mu |b_1|^\kappa = \mathbb{E}_\mu(V_1^{\kappa/2}) \mathbb{E}_\mu |\xi_1|^\kappa$, and $\xi_1$ is a standard Gaussian variable and $V_1$ is bounded.

*Second term.* We have

$$\int_0^{e^t} u^\kappa \mathbb{P}_\mu(u - u^\varepsilon < \Phi_1 R_0 \le u, X_\delta = i)\, du$$

$$= \int_0^{e^t} u^\kappa \mathbb{P}_\mu(\Phi_1 R_0 > u - u^\varepsilon, X_\delta = i)\, du$$

$$\quad - \int_0^{e^t - e^{t\varepsilon}} u^\kappa \mathbb{P}_\mu(\Phi_1 R_0 > u, X_\delta = i)\, du$$

$$\le \int_0^{e^t} u^\kappa [1 - \mathbf{1}_{u \ge 1}(u - u^\varepsilon)^\kappa (1 - \varepsilon u^{\varepsilon - 1})] \mathbb{P}_\mu(\Phi_1 R_0 > u - u^\varepsilon, X_\delta = i)\, du.$$

Set $0 < r < \kappa$. As $\Phi_1$ is bounded, there exists a positive constant $c$ such that for all $u > 0$ we have

$$\mathbb{P}_\mu(\Phi_1 R_0 > u, X_1 = i) \le c \frac{\mathbb{E}_\mu |R_0|^r}{u^r},$$

which is bounded by Proposition 6. Thus we get

$$(15) \qquad \int_0^{e^t} u^\kappa \mathbb{P}_\mu(u - u^\varepsilon < \Phi_1 R_0 \le u, X_1 = i)\, du \le C e^{t(\kappa - r + \varepsilon - 1)},$$

where $C$ is a positive constant. Now set $\beta = \max\{\kappa + \varepsilon - r; 1 + \kappa - \kappa\varepsilon\} \in\, ]0; 1[$. Then (13)–(15) yield $e^t G_i^1(t) \le c e^{t\beta}$ for all $t > 0$. Thus $G_i^1(t) \le c e^{t(\beta - 1)}$ is directly Riemann integrable on $\mathbb{R}_+$.

6.5. *Tail of the distribution.* We have now proved that $F$ and $G$ satisfy the assumptions of Theorem A. Thus we get, for all $i, t$,

$$(16) \qquad Z_i(t) \underset{t \to \infty}{\longrightarrow} c m_i \sum_{j=1}^N u_j \int_{-\infty}^\infty G_j(y)\, dy.$$

Summing up these terms, we get

$$(17) \qquad z(t) \underset{t \to \infty}{\longrightarrow} c \sum_{j=1}^N u_j \int_{-\infty}^\infty G_j(y)\, dy,$$

as $\sum m_i = 1$. We still have to prove that this limit is nonzero.

**7. Proof of Theorem 2, part II.** Now we are going to prove that there exists a positive constant $C$ such that $t^\kappa \mathbb{P}_\mu(|R_1| > t) \geq C > 0$ when $t$ tends to infinity. First, we give a lower bound of this probability involving the products $\Phi_1 \cdots \Phi_n$, and then we study the asymptotic behavior of such products.

7.1. *Lower bound for $\nu\{x : |x| > t\}$.* The following proof is adapted from Goldie (1991) and de Saporta (2004).

PROPOSITION 8. *There exist $\varepsilon > 0$ and a corresponding positive constant $C$ such that for large enough $t$ we have*

$$\mathbb{P}_\mu(|R_1| > t) \geq C \mathbb{P}_\mu\left( \sup_n (\Phi_1 \cdots \Phi_n) > \frac{2t}{\varepsilon} \right).$$

For the i.i.d. case, the key to such a lower bound is an inequality established in Grincevičius (1980) that extends Lévy's symmetrization inequality [see Chow and Teicher (1978)]. Here we need first to extend this inequality.

Recall that $R_1 = \sum_{k=0}^\infty \Phi_1 \Phi_0 \cdots \Phi_{2-k} b_{1-k}$. For all $n \geq 1$, we set

$$R_1^n = \sum_{k=0}^{n-1} \Phi_1 \Phi_0 \cdots \Phi_{2-k} b_{1-k} \quad \text{and} \quad \Pi_n = \Phi_1 \Phi_0 \cdots \Phi_{2-n}.$$

If $x$ is a $\sigma(X_t, W_t, a \leq t \leq b)$-measurable random variable, let $\operatorname{med}_i(x)$ be a median of $x$ conditionally to $X_b = i$ and $\operatorname{med}_-(x) = \min_i\{\operatorname{med}_i(x)\}$.

LEMMA 3. *For all $t > 0$ and $n \geq 1$, we have*

$$\mathbb{P}_\mu\left( \max_{1 \leq j \leq n} \left\{ R_1^j + \Pi_j \operatorname{med}_-\left( \frac{R_1^n - R_1^j}{\Pi_j} \right) \right\} > t \right) \leq 2\mathbb{P}_\mu(R_1^n > t).$$

PROOF. Set $T = \inf\{j \leq n \text{ t.q. } R_1^j + \Pi_j \operatorname{med}_-(\Pi_j^{-1}(R_1^n - R_1^j)) > t\}$ if this set is not empty, $n + 1$ otherwise, and $B_j = \{\operatorname{med}_-(\Pi_j^{-1}(R_1^n - R_1^j)) \leq \Pi_j^{-1}(R_1^n - R_1^j)\}$. The event $(T = j)$ is in the $\sigma$-field generated by $(X_t, W_t, (1 - j) \leq t \leq 1)$, and $B_j$ is in the $\sigma$-field generated by $(X_t, W_t, (1 - n) \leq t \leq (1-j))$. Therefore these events are independent conditionally to $X_{(1-j)}$. Besides, for all $i$ and $j$ we have $\mathbb{P}_\mu(B_j | X_{(1-j)} = i) \geq \mathbb{P}_\mu(\operatorname{med}_i(\Pi_j^{-1}(R_1^n - R_1^j)) \leq \Pi_j^{-1}(R_1^n - R_1^j) | X_{(1-j)} = i) \geq 1/2$. Thus, as the products $\Pi_j$ are positive or zero, we have

$$\mathbb{P}_\mu(R_1^n > t) \geq \mathbb{P}_\mu\left( \bigcup_{j=1}^n [B_j \cap (T = j)] \right)$$

$$= \sum_{j=1}^n \sum_{i=1}^N \mathbb{P}_\mu(B_j | X_{(1-j)} = i) \mathbb{P}(T = j | X_{(1-j)} = i) \mu(i)$$

$$\geq \frac{1}{2} \mathbb{P}_\mu(T \leq n)$$

$$= \frac{1}{2} \mathbb{P}_\mu \left( \max_{1 \leq j \leq n} \left\{ R_1^j + \Pi_j \operatorname{med} \left( \frac{R_1^n - R_1^j}{\Pi_j} \right) \right\} > t \right). \qquad \square$$

Under our assumptions, $R_1^n$ tends to $R_1$ when $n$ tends to infinity, and for fixed $j$, $\Pi_j^{-1}(R_1^n - R_1^j)$ converges to a random variable $\widehat{R}$ that has the same distribution as $R_1$. Set $m_0 = \operatorname{med}_-(R_1) = \min_i \{ \operatorname{med}(R_1|X_1 = i) \} = \operatorname{med}_-(\widehat{R})$, and letting $n$ tend to infinity in Lemma 3, we get, for all $t > 0$,

$$\mathbb{P}_\mu \left( \sup_j \{ R_1^j + \Pi_j m_0 \} > t \right) \leq 2 \mathbb{P}_\mu(R_1 > t).$$

Replacing $R_1$ by $-R_1$ yields a similar formula; thus, for all $t > 0$ we get

(18) $$\mathbb{P}_\mu \left( \sup_j |R_1^j + \Pi_j m_0| > t \right) \leq 2 \mathbb{P}_\mu(|R_1| > t).$$

Furthermore, as proved in Goldie [(1991), page 157], for all $t > |m_0|$ we have

$$\mathbb{P}_\mu \left( \sup_n \{ R_1^n + \Pi_n m_0 \} > t \right)$$

$$\geq \mathbb{P}_\mu(\exists n \text{ s.t. } |(R_1^{n+1} + \Pi_{n+1} m_0) - (R_1^n + \Pi_n m_0)| > 2t),$$

where $R_1^0 = 0$ and $\Pi_0 = 1$. But we have

$$(R_1^{n+1} + \Pi_{n+1} m_0) - (R_1^n + \Pi_n m_0)$$

$$= \Phi_1 \Phi_0 \cdots \Phi_{2-n} b_{1-n} + (\Pi_{n+1} - \Pi_n) m_0$$

$$= \Pi_n (b_{1-n} + (\Phi_{1-n} - 1) m_0).$$

Thus (18) yields, for all $\varepsilon > 0$,

$$\mathbb{P}_\mu(|R_1| > t) \geq \frac{1}{2} \mathbb{P}_\mu(\exists n \text{ s.t. } |\Pi_n(b_{1-n} + (\Phi_{1-n} - 1) m_0)| > 2t)$$

(19)

$$\geq \frac{1}{2} \mathbb{P}_\mu \left( \exists n \text{ s.t. } |\Pi_n| > \frac{2t}{\varepsilon} \text{ and } |b_{1-n} + (\Phi_{1-n} - 1) m_0| > \varepsilon \right).$$

Now we give an extension of Feller–Chung's inequality adapted to the present case [see Chow and Teicher (1978)]:

LEMMA 4. *For all $t > |m_0|$ and $\varepsilon > 0$, we have*

$$\mathbb{P}_\mu \left( \exists n \text{ s.t. } |\Pi_n| > \frac{2t}{\varepsilon} \text{ and } |b_{1-n} + (\Phi_{1-n} - 1) m_0| > \varepsilon \right)$$

$$\geq \min_{1 \leq i \leq N} \mathbb{P}_i(|b_0 + (\Phi_0 - 1) m_0| > \varepsilon) \mathbb{P}_\mu \left( \exists n \text{ s.t. } |\Pi_n| > \frac{2t}{\varepsilon} \right).$$

PROOF.   Set $A_0 = \varnothing$, $A_n = \{|\Pi_n| > 2t\varepsilon^{-1}\}$ and $B_n = \{|b_{1-n} + (\Phi_{1-n} - 1)m_0| > \varepsilon\}$. Conditionally to $X_{(1-n)}$, $B_n$ is independent of $A_0, \ldots, A_n$. Thus we have

$$\mathbb{P}_\mu\left(\bigcup_{n=1}^\infty [A_n \cap B_n]\right)$$

$$= \sum_{n=1}^\infty \mathbb{P}_\mu\left(B_n \cap A_n \bigcap_{j=0}^{n-1} [B_j \cap A_j]^c\right)$$

$$\geq \sum_{n=1}^\infty \mathbb{P}_\mu\left(B_n \cap A_n \bigcap_{j=0}^{n-1} A_j^c\right)$$

$$= \sum_{n=1}^\infty \sum_{i=1}^N \left[\mathbb{P}_\mu(B_n|X_{(1-n)} = i)\mathbb{P}_\mu\left(A_n \bigcap_{j=0}^{n-1} A_j^c \Big| X_{(1-n)} = i\right)\mu(i)\right],$$

where $A^c$ denotes the complementary set of $A$. But, by stationarity we have $\mathbb{P}_\mu(B_n|X_{(1-n)} = i) = \mathbb{P}_i(|b_0 + (\Phi_0 - 1)m_0| > \varepsilon)$. Thus we get

$$\mathbb{P}_\mu\left(\bigcup_{n=1}^\infty [A_n \cap B_n]\right) \geq \min_{1 \leq i \leq N} \mathbb{P}_i(|b_0 + (\Phi_0 - 1)m_0| > \varepsilon)\mathbb{P}_\mu\left(\bigcup_{n=1}^\infty A_n\right). \quad \square$$

Now we can give the proof of Proposition 8.

PROOF OF PROPOSITION 8.   Equation (19) and Lemma 4 yield, for all $t > |m_0|$ and $\varepsilon > 0$,

$$\mathbb{P}_\mu(|R_1| > t) \geq \frac{1}{2} \min_{1 \leq i \leq N} \mathbb{P}_i(|b_0 + (\Phi_0 - 1)m_0| > \varepsilon)\mathbb{P}_\mu\left(\exists n \text{ s.t. } |\Pi_{n-1}| > \frac{2t}{\varepsilon}\right).$$

We have $b_0 = V_0^{1/2}\xi_0$, $V_0$ and $\Phi_0$ are bounded, but $\xi$ is not bounded as it is a Gaussian variable. Thus equality $b_0 + (\Phi_0 - 1)m_0 = 0$ cannot hold $\mathbb{P}_i$-almost surely. Thus we can find $\varepsilon > 0$ such that $\min_{1 \leq i \leq N} \mathbb{P}_i(|b_0 + (\Phi_0 - 1)m_0| > \varepsilon) > 0$. Hence, as expected there is a constant $C > 0$ such that for all $t > |m_0|$, we have

$$\mathbb{P}_\mu(|R_1| > t) \geq C\mathbb{P}_\mu\left(\sup_n |\Pi_n| > \frac{2t}{\varepsilon}\right). \quad \square$$

7.2. *Asymptotic behavior of the products* $\Phi_1 \cdots \Phi_n$.   To estimate the probability $\mathbb{P}_\mu(\sup_n |\Pi_n| > t)$, we use the ladder height method given by Feller (1966) for the study of the maximum of random walks.

7.2.1. *Notation.* First we introduce some notation. Set $S_0 = 0$ and for all positive $n$, we set

$$S_n = \sum_{k=1}^{n} \log(\Phi_{2-k}) = \log \Pi_n = \int_{(1-n)}^{1} a(X_u)\, du.$$

The *first ladder epoch* of this random walk is $\tau = \tau_1 = \inf\{n \geq 1 \text{ s.t. } S_n > 0\}$, and the *first ladder height* is $S_\tau$. We denote by $H(t)$ the matrix of distributions of $S_\tau$ with the following coordinates:

$$H_{ij}(t) = \mathbb{P}_\mu(\tau < \infty, S_\tau \leq t, X_{(1-\tau)} = j | X_1 = i).$$

As $S_\tau > 0$, $H$ is distributed on the positive half-line. Moreover, $S_\tau > 0$, $S_{1-\tau} \leq 0$ and the $\Phi_n$ are bounded, thus we have $S_\tau \leq \sup \log \Phi_n \leq \sup a(i) < \infty$, and $H$ has bounded support.

We define also the $n$th ladder epoch by $\tau_n = \inf\{n > \tau_{n-1} \text{ s.t. } S_n > S_{\tau_{n-1}}\}$, and $S_{\tau_n}$ is the corresponding ladder height. We check that we have

$$H_{ij}^{(n)}(t) = \mathbb{P}_\mu(\tau_n < \infty, S_{\tau_n} \leq t, X_{(1-\tau_n)} = j | X_1 = i),$$

where $H^{(n)}$ is the $n$-fold convolution of $H$. Let $\Psi = \sum_{n=0}^{\infty} H^{(n)}$ be the renewal function associated with $H$.

7.2.2. *The random walk $S_{\tau_n}$.* To investigate the asymptotic behavior of $(S_{\tau_n})$ we are going to use a renewal theorem as in Feller (1966) for the i.i.d. case, namely, Theorem B. We want to apply it for $F = H$ and $\alpha = \kappa$, thus we have to prove that $H$ satisfies its assumptions.

As $H(0) = 0$, we have $\rho(H(0)) < 1$, thus the first assumption is true. In addition, $H_{ij}$ are probability measures, therefore $H$ is finite. $H$ has bounded support because $S_{\tau-1} \leq 0$, $S_\tau > 0$ and $\Phi$ is bounded. Thus $\widehat{B}$, the expectation of $H_\kappa(\infty) = \int_0^\infty e^{-\kappa u} H(du)$ is well defined. Proposition 7 yields again that $H$ is also nonlattice.

*Irreducibility and aperiodicity.* For all $i, j$ in $E$, we have

$$H_{ij}(\infty) = \mathbb{P}_\mu(\tau < \infty, \ X_{1-\tau} = j | X_1 = i)$$

$$\geq \mathbb{P}_\mu(\tau = 1, \ X_0 = j | X_1 = i)$$

$$= \mathbb{P}_j(\log \Phi_1 > 0, \ X_1 = i) \frac{\mu(j)}{\mu(i)}$$

$$= \mathbb{P}_j\left(\int_0^1 a(X_u)\, du > 0, \ X_1 = i\right) \frac{\mu(j)}{\mu(i)},$$

and Proposition 7 implies that the last term is positive as $0 \in\, ]a_m; a_M[$. Thus the second assumption of Theorem B is valid. We have also proved that for all $i$ and $j$ we have $0 = H_{ij}(0) < H_{ij}(\infty)$, so that the third assumption is also valid.

*Spectral radius of $H_\kappa(\infty)$.* Now we define a new probability law $\mathbb{P}_\kappa$ on $\Omega \times \Theta$. For all bounded $\mathcal{A} \times \mathcal{B}$-measurable functions $f$ which first coordinate depends only on $(X_t, (1-n) \leq t \leq 1)$, we set

$$\mathbb{P}_\kappa(f) = \mathbb{E}_\kappa(f)$$
$$= \frac{\mathbb{E}_\mu(f(\Phi_1, \ldots, \Phi_{2-n}, \theta)(\Phi_1 \cdots \Phi_{2-n})^\kappa)}{\mathbb{E}_\mu((\Phi_1 \cdots \Phi_{2-n})^\kappa)}.$$

Set $H_\kappa(t) = \int_0^t e^{-\kappa u} H(du)$. We have

$$(H_\kappa)_{ij}(t) = \frac{\mathbb{P}_\kappa(\tau < \infty, S_\tau \leq t, X_{(1-\tau)} = j | X_1 = i)}{\mathbb{E}_\mu((\Phi_1 \cdots \Phi_{1-\tau})^\kappa, \tau < \infty)}$$
$$= \frac{(H_\kappa^\sharp)_{ij}(t)}{\mathbb{E}_\mu((\Phi_1 \cdots \Phi_{1-\tau})^\kappa, \tau < \infty)},$$

where $(H_\kappa^\sharp)_{ij}(t) = \mathbb{P}_\kappa(\tau < \infty, S_\tau \leq t, X_{(1-\tau)} = j | X_1 = i)$ describes the behavior of the ladder heights of our random walk under the new probability law $\mathbb{P}_\kappa$.

The computation we made in the proof of Proposition 3 yields

$$\frac{\partial}{\partial r}\Big|_{r=\kappa} \log(\rho(A_r)) = \lim_{n \to \infty} \frac{1}{n} \mathbb{E}_\kappa\left(\sum_{i=1}^n \log \Phi_i\right)$$
$$= \mathbb{E}_\kappa(\log \Phi_1).$$

But we have $\log \rho(A_0) = \log \rho(A_\kappa) = 0$; this function is convex (Corollary 2) and its right-hand derivative at 0 is negative (Proposition 3). Thus its left-hand derivative at $\kappa$ is positive, that is, $\mathbb{E}_\kappa(\log \Phi_1) > 0$. Under the law $\mathbb{P}_\kappa$ our random walk thus drifts to $+\infty$, hence for all $n$ and $i$, we have $(\mathbb{P}_\kappa)_i(\tau_n < \infty) = 1$ and $H^\sharp$ is a stochastic matrix, therefore its spectral radius equals to 1.

For all $n$, we have

$$H_\kappa^{(n)}(\infty) = (H_\kappa(\infty))^n = \frac{(H_\kappa^\sharp(\infty))^n}{\mathbb{E}_\mu((\Phi_1 \cdots \Phi_{2-\tau_n})^\kappa, \tau < \infty)},$$

thus $\rho(H_\kappa(\infty)) = \lim(\mathbb{E}_\mu((\Phi_1 \cdots \Phi_{2-\tau_n})^\kappa, \tau < \infty))^{-1/n}$ and we now have to prove that this limit equals to 1. But for all $n$, we have $\tau_n \geq n$, and the event $(\tau_n = k)$ depends only on $(X_t, (1-k) \leq t \leq 1)$. Thus we have

$$\mathbb{E}_\mu((\Phi_1 \cdots \Phi_{1-\tau_n})^\kappa, \tau_n < \infty)$$

(20)
$$= \sum_{k=n}^\infty \mathbb{E}_\mu((\Phi_1 \cdots \Phi_{1-k})^\kappa, \tau_n = k)$$

$$= \sum_{k=n}^\infty \mathbb{P}_\kappa(\tau_n = k)\mathbb{E}_\mu((\Phi_1 \cdots \Phi_{1-k})^\kappa).$$

Set $\varepsilon > 0$. For large enough $n$, our choice of $\kappa$ and (5) and (6) yield

$$\mu A_\kappa^n \mathbf{1} - \varepsilon \leq \mathbb{E}_\mu((\Phi_1 \cdots \Phi_{1-n})^\kappa) \leq \mu A_\kappa^n \mathbf{1} + \varepsilon.$$

Thus for large enough $n$, (20) yields

$$(\mu A_\kappa^n \mathbf{1} - \varepsilon) \sum_{k=n}^\infty \mathbb{P}_\kappa(\tau_n = k)$$

$$\leq \mathbb{E}_\mu((\Phi_1 \cdots \Phi_{1-\tau_n})^\kappa, \tau_n < \infty) \leq (\mu A_\kappa^n \mathbf{1} + \varepsilon) \sum_{k=n}^\infty \mathbb{P}_\kappa(\tau_n = k),$$

and as $\mathbb{P}_\kappa(\tau_n < \infty) = 1$, we have

$$\mu A_\kappa^n \mathbf{1} - \varepsilon \leq \mathbb{E}_\mu((\Phi_1 \cdots \Phi_{1-\tau_n})^\kappa, \tau_n < \infty) \leq \mu A_\kappa^n \mathbf{1} + \varepsilon.$$

Thus as $n \to \infty$ we have, with the notation of Corollary 1, $\mathbb{E}_\mu(\Phi_1 \cdots \Phi_{1-\tau_n})^\kappa \sim \mu B_\kappa \mathbf{1}$. Hence we have, as expected, $\mathbb{E}_\mu((\Phi_1 \cdots \Phi_{1-\tau_n})^\kappa, \tau_n < \infty)^{1/n} \to 1$.

Thus all the assumptions of Theorem B are valid here. We are going to use it in the following part.

7.2.3. *Asymptotic behavior of the maximum.* Let $M = \sup_n S_n = \sup_n S_{\tau_n}$ be the maximum of our random walk. Using the definition of $H$, we get, for all $1 \leq i \leq N$,

$$\mathbb{P}_\mu(M \leq t | X_1 = i)$$

$$= \sum_{n=1}^\infty \mathbb{P}_\mu(\tau_n < \infty, S_{\tau_n} \leq t, \tau_{n+1} = \infty | X_1 = i)$$

$$= \sum_{n=1}^\infty \sum_{j=1}^N \mathbb{P}_\mu(\tau_n < \infty, S_{\tau_n} \leq t, \tau_{n+1} = \infty, X_1 = i | X_{(1-\tau_n)} = j) \frac{\mu(j)}{\mu(i)}$$

$$(21) \qquad = \sum_{n=1}^\infty \sum_{j=1}^N [\mathbb{P}_\mu(\tau_n < \infty, S_{\tau_n} \leq t, X_{(1-\tau_n)} = j | X_1 = i)$$

$$\times (1 - \mathbb{P}_\mu(\tau_{n+1} < \infty | X_{(1-\tau_n)} = j))]$$

$$= \sum_{n=1}^\infty \sum_{j=1}^N \left[ H_{ij}^{(n)}(t) \left( 1 - \sum_{k=1}^N H_{jk}(\infty) \right) \right]$$

$$= \sum_{j=1}^N \left[ \Psi_{ij}(t) \left( 1 - \sum_{k=1}^N H_{jk}(\infty) \right) \right].$$

Theorem B applied to (21) yields, when $t$ tends to infinity,

$$1 - \mathbb{P}_\mu(M \leq t | X_1 = i)$$

$$(22) \quad \begin{aligned}
&= \sum_{j=1}^{N} \left[ \left( 1 - \sum_{k=1}^{N} H_{jk}(\infty) \right) \int_{t}^{\infty} e^{-\kappa u} (e^{\kappa u} \Psi_{ij})(du) \right] \\
&\overset{t \to \infty}{\sim} \sum_{j=1}^{N} \left[ \left( 1 - \sum_{k=1}^{N} H_{jk}(\infty) \right) \int_{t}^{\infty} e^{-\kappa u} \widehat{c}\, \widehat{m}_i \widehat{u}_j \, du \right] \\
&= \sum_{j=1}^{N} \left[ \left( 1 - \sum_{k=1}^{N} H_{jk}(\infty) \right) \widehat{c}\, \widehat{m}_i \widehat{u}_j \right] e^{-\kappa t},
\end{aligned}$$

where $\widehat{m}$ and $\widehat{u}$ are right and left eigenvectors of $H_\kappa(\infty)$ with positive coordinates with the same normalization as in Section 4, and $\widehat{c} = ({}^t \widehat{u} \widehat{B} \widehat{m})^{-1} > 0$.

7.3. *Conclusion.* We still have to prove that there is a $j \leq N$ such that $1 - \sum_{k=1}^{N} H_{jk}(\infty) > 0$. But the mapping $r \longmapsto H_r(\infty) = \int_0^\infty e^{ru} H(du)$ is clearly increasing component-wise. As these matrices are nonnegative and irreducible, Corollaries 8.1.19 and 8.1.20 of Horn and Johnson (1985) imply that the mapping $r \longmapsto \rho(H_r(\infty))$ is also increasing. As $\rho(H_\kappa(\infty)) = 1$, we have $\rho(H_0(\infty)) = \rho(H(\infty)) < 1$. This is a substochastic, nonstochastic matrix, thus there exists a $j$ such that we have $1 - \sum_{k=1}^{N} H_{jk}(\infty) > 0$.

We have now proved that the right-hand side term in (22) is positive, thus there is a constant $C > 0$ such that, when $t$ tends to infinity, we have

$$(23) \quad e^{\kappa t} \mathbb{P}_\mu(M > t) \geq \sum_{i=1}^{N} e^{\kappa t} \mathbb{P}_\mu(M > t | X_1 = i) \mu(i) \geq C.$$

Putting together this result and Proposition 8, we get, for large enough $t$,

$$(24) \quad t^\kappa \mathbb{P}_\mu(|R_1| > t) \geq K > 0.$$

With the notation of Theorem 2, it means that $L > 0$, which ends the proof of this theorem.

8. **Determination of $\kappa$.** Set $s_1 = \min\{\lambda(i)a(i)^{-1} | a(i) > 0\}$, and let $M_s$ be the matrix with components $\{q(i,j)\lambda(i)(\lambda(i) - sa(i))^{-1}\}$. This matrix is well defined for all $s < s_1$. We can precisely compare the spectral radius of $A_s$ and that of $M_s$.

PROPOSITION 9. *For all $0 < s < s_1$, we have $\rho(M_s) < 1$ if and only if $\rho(A_s) < 1$, and we have $\rho(M_s) > 1$ if and only if $\rho(A_s) > 1$.*

PROOF. Suppose that $\rho(M_s) < 1$. $M_s$ is a positive irreducible matrix as $q$ is, $\lambda$ being positive and $s < s_1$. Thus the Perron–Frobenius theorem [see, e.g., Horn and Johnson (1985)] gives the existence of a vector $\varphi$ with positive coordinates such that $M_s\varphi = \rho(M_s)\varphi < \varphi$. Hence for all $i$ in $E$, we have

$$\varphi(i) > \sum_j q(i,j)\frac{\lambda(i)}{\lambda(i) - sa(i)}\varphi(j),$$

that we can rewrite, since $s < s_1$, as

$$(25) \qquad (sa(i) - \lambda(i))\varphi(i) + \lambda(i)\sum_j q(i,j)\varphi(j) < 0.$$

Proposition 4 enables us to choose a small enough $\delta$ such that (8) is valid here. Equation (25) thus yields

$$A_s\varphi(i) = [1 + \delta(sa(i) - \lambda(i))]\varphi(i) + \delta\lambda(i)\sum_{j\neq i}[q(i,j)\varphi(j)] + o(\delta)$$

$$= \varphi(i) + \delta\left[(sa(i) - \lambda(i))\varphi(i) + \lambda(i)\sum_j q(i,j)\varphi(j)\right] + o(\delta)$$

$$< \varphi(i).$$

Thus component-wise we get $A_s\varphi < \varphi$, which implies that $\rho(A_s) < 1$. The proof that $\rho(M_s) > 1$ implies $\rho(A_s) > 1$ runs the same, the inequalities being reversed.

Suppose now that $\rho(A_s) < 1$. $A_s$ is a positive irreducible matrix, thus the Perron–Frobenius theorem gives the existence of a vector $\psi$ with positive coordinates such that $A_s\psi = \rho(A_s)\psi < \psi$. Hence for all $i$ in $E$, and small enough $\delta$, we have

$$\delta\left[(sa(i) - \lambda(i))\psi(i) + \lambda(i)\sum_j q(i,j)\psi(j)\right] + o(\delta) = A_s\psi(i) - \psi_i$$

$$< 0.$$

Hence, for all $i$, we get $(sa(i) - \lambda(i))\psi(i) + \lambda(i)\sum_j q(i,j)\psi(j) < 0$, or, as $s < s_1$,

$$\psi(i) > \frac{\lambda(i)}{\lambda(i) - sa(i)}\sum_j q(i,j)\psi(j),$$

and thus $M_s\psi < \psi$. As $M_s$ is a positive matrix, we conclude that $\rho(M_s) < 1$. Here again the proof that $\rho(A_s) > 1$ implies $\rho(M_s) > 1$ runs the same with reversed inequalities. $\quad\square$

PROPOSITION 10. *The spectral radius of $M_s$ tends to infinity when $s$ tends to $s_1$.*

Proof. Set $i_0 \in E$ such that $\lambda(i_0)a(i_0)^{-1} = s_1$, and $e_{i_0}$ the row vector with zero coordinates except the $i_0$th which is set to be 1. Set $v_{i_0} = \lambda(i_0)(\lambda(i_0) - sa(i_0))^{-1}$. We have $e_{i_0}M_s = v_{i_0}q(i_0, \cdot) \geq v_{i_0}e_{i_0}$ as $q$ is a positive matrix. As $M_s$ is also positive, for all $s < s_1$, we get $\rho(M_s) \geq v_{i_0} = \lambda(i_0)(\lambda(i_0) - sa(i_0))^{-1}$. Hence this spectral radius tends to infinity when $s$ tends to $s_1$.  □

Corollary 4. *There is a unique $s \in\, ]0; s_1[$ such that $\rho(M_s) = 1$, and this $s$ equals the unique $\kappa$ such that $\rho(A_\kappa) = 1$.*

Proof. For all $s < \kappa$, we have $\rho(A_s) < 1$ by Corollary 3; thus Proposition 9 yields $\rho(M_s) < 1$ for all $0 < s < \min\{\kappa, s_1\}$. As $\rho(M_s) \to \infty$ as $s$ tends to $s_1$, we also have $\rho(A_s) > 1$ for $s$ close to $s_1$. Therefore $\kappa < s_1$, and $\rho(A_s) > 1$ for all $\kappa < s < s_1$. Hence $\rho(M_s) > 1$ for all $\kappa < s < s_1$. As $M_s$ has continuous coordinates, its spectral radius is also continuous; thus $\rho(M_\kappa) = 1$ and $\kappa$ is the only value of $s \in\, ]0; s_1[$ satisfying this equation.  □

We now give an illustration by computing the value of $\kappa$ when $E = \{1, 2\}$. The jump kernel $q$ then equals to

$$q = \begin{pmatrix} 0 & 1 \\ 1 & 0 \end{pmatrix},$$

and the invariant law of the process $X$ is $\mu = (\lambda(2), \lambda(1))/(\lambda(1) + \lambda(2))$. We suppose that $a(1)$ or $a(2)$ is positive. Condition 2 becomes

(26) $$\lambda(1)a(2) + \lambda(2)a(1) < 0.$$

For all $i$ in $E$, set $r_i = \frac{a(i)}{\lambda(i)}$. We have $r_1 + r_2 < 0$, $r_1 r_2 > 0$ and $s_1 = \max\{r_1^{-1}, r_2^{-1}\}$. For $s \in [0, s_1[$, the matrix $M_s$ equals to

$$M_s = \begin{pmatrix} 0 & \dfrac{1}{1 - sr_1} \\ \dfrac{1}{1 - sr_2} & 0 \end{pmatrix},$$

and its spectral radius is $[(1 - sr_1)(1 - sr_2)]^{-1/2}$. It equals to 1 for $\kappa = r_1^{-1} + r_2^{-1} = \lambda(2)a(2)^{-1} + \lambda(1)a(1)^{-1}$.

## REFERENCES

Athreya, K. and Rama Murthy, K. (1976). Feller's renewal theorem for systems of renewal equations. *J. Indian Inst. Sci.* **58** 437–459. MR436366

Barndorff-Nielsen, O. E. and Shephard, N. (2001). Non-Gaussian Ornstein–Uhlenbeck-based models and some of their uses in financial economics. *J. Roy. Statist. Soc. Ser. B* **63** 167–241. MR1841412

Basak, G., Bisi, A. and Ghosh, M. K. (1996). Stability of random diffusion with linear drift. *J. Math. Anal. Appl.* **202** 604–622. MR1406250

Brockwell, P. J. (2001). Lévy-driven CARMA processes. *Ann. Inst. Statist. Math.* **53** 113–124. MR1820952

Chow, S. C. and Teicher, H. (1978). *Probability Theory. Independence, Interchangeability, Martingales.* Springer, New York. MR513230

Crump, K. (1970). On systems of renewal equations. *J. Math. Anal. Appl.* **30** 425–434. MR257678

de Saporta, B. (2003). Renewal theorem for a system of renewal equations. *Ann. Inst. H. Poincaré Probab. Statist.* **39** 823–838. MR1997214

de Saporta, B. (2004). Tail of the stationary solution of the stochastic equation $y_{n+1} = a_n y_n + b_n$ with Markovian coefficients. *Stochastic Process. Appl.* To appear.

Feller, W. (1966). *Introduction to Probability Theory and Its Applications* **2**. Wiley, New York. MR210154

Goldie, C. (1991). Implicit renewal theory and tails of solutions of random equations. *Ann. Appl. Probab.* **1** 26–166. MR1097468

Grincevičius, A. K. (1980). Products of random affine transformations. *Lithuanian Math. J.* **20** 279–282. MR605960

Guyon, X., Iovleff, S. and Yao, J.-F. (2004). Linear diffusion with stationary switching regime. *ESAIM Probab. Statist.* **8** 25–35. MR2085603

Hamilton, J. (1989). A new approach to the economic analysis of nonstationary time series and the business cycle. *Econometrica* **57** 151–173. MR996941

Hamilton, J. (1990). Analysis of time series subject to changes in regime. *J. Econometrics* **45** 39–70. MR1067230

Horn, R. and Johnson, C. (1985). *Matrix Analysis.* Cambridge Univ. Press. MR832183

Karatzas, I. and Shreve, S. (1991). *Brownian Motion and Stochastic Calculus.* Springer, New York. MR1121940

Kesten, H. (1973). Random difference equations and renewal theory for products of random matrices. *Acta Math.* **131** 207–248. MR440724

Le Page, E. (1983). Théorèmes de renouvellement pour les produits de matrices aléatoires. Equations aux différences aléatoires. In *Séminaires de Probabilités de Rennes* 116. Publ. Sem. Math., Univ. Rennes 1, Rennes. MR863321

Norris, J. (1998). *Markov Chains.* Cambridge Univ. Press. MR1600720

IRMAR
Université de Rennes 1
Campus de Beaulieu
35042 Rennes Cedex
France
e-mail: benoite.de-saporta@math.univ-rennes1.fr
e-mail: jian-feng.yao@univ-rennes1.fr